\documentclass[11pt]{amsart}
\usepackage{amssymb,latexsym,amsmath,setspace}
\input{xypic}
\xyoption{all}
\numberwithin{equation}{section}

\newcommand{\Q}{{\mathbb Q}}
\newcommand{\R}{{\mathbb R}}
\newcommand{\Z}{{\mathbb Z}}
\newcommand{\C}{{\mathbb C}}

\newcommand{\A}{{\mathbb A}}
\newcommand{\I}{{\mathbb I}}

\newcommand{\g}{{\mathfrak g}}

\renewcommand{\O}{{\mathcal O}}

\newcommand{\M}{{\mathcal M}}

\newcommand{\GL}{{\rm GL}}
\newcommand{\SO}{{\rm SO}}
\newcommand{\Sp}{{\rm Sp}}
\newcommand{\SL}{{\rm SL}}
\newcommand{\PU}{{\rm PU}}

\def\cE{\mathcal{E}}

\def\Orth{{\rm O}}
\def\Unit{{\rm U}}
\def\gl{\mathfrak{gl}}
\def\k{\mathfrak k}
\def\tM{\widetilde{\mathcal{M}}}
\def\Coh{{\rm Coh}} 
\def\autc{{\rm Aut}(\C)}
\def\Ind{{\rm Ind}} 
\def\cG{{\mathcal G}}

\numberwithin{equation}{subsection}
\newtheorem{theorem}[equation]{Theorem}
\newtheorem{cor}[equation]{Corollary}

\newtheorem{prop}[equation]{Proposition}
\newtheorem{con}[equation]{Conjecture}

\newtheorem{defn}[equation]{Definition}

\newtheorem{prob}[equation]{Problem}

\begin{document}

\title{Endoscopy and the Cohomology of GL(n)}
\date{\today}
\subjclass[2010]{Primary: 11F75; Secondary: 11F70, 22E55}
\author{\bf Chandrasheel Bhagwat \ \ \& \ \ A. Raghuram}
\address{Indian Institute of Science Education and Research, Dr.\,Homi Bhabha Road, Pashan, Pune 411008,  INDIA.}
\email{cbhagwat@iiserpune.ac.in, \ raghuram@iiserpune.ac.in}

\maketitle

Let $F$ be a totally real field. In Grobner--Raghuram~\cite{Grobner-Raghuram} the special values of the standard $L$-function of a cuspidal automorphic representation $\pi$ of $\GL(2n)/F$ of cohomological type and admitting a Shalika model were studied. As explained via examples in {\it loc.\,cit.,} the geometric condition of being cohomological and the analytic condition of admitting a Shalika model are entirely independent of each other, and there is no {\it a priori} reason why such a $\pi$ should even exist. The purpose of this article is to address such existence questions; indeed, we prove that such a $\pi$ exists.  More generally,  let $G = {\rm Res}_{F/\Q}(\GL_n)$ where $F$ now is any number field, and let $S^G_{K_f}$ denote an ad\`elic locally symmetric space for some level structure $K_f.$ Let  $\M_{\mu,\C}$ be an algebraic irreducible representation of $G(\R)$ and we let $\tM_{\mu,\C}$ denote the associated sheaf on $S^G_{K_f}.$ The basic problem addressed in this paper is to classify the data $(F,n,\mu)$ for which cuspidal cohomology of $G$ with $\mu$-coefficients, denoted $H^{\bullet}_{\rm cusp}(S^G_{K_f}, \tM_{\mu,\C})$, is nonzero for some $K_f.$ We 
prove nonvanishing of cuspidal cohomology when $F$ is a totally real field or a totally imaginary quadratic extension of a totally real field, and also for a general number field but when $\mu$ is a parallel weight. The proof in the totally real case involves Arthur's endoscopic classification \cite{Arthur} of the discrete spectrum for classical groups, together with Clozel's results on globalizing discrete series representations via limit multiplicity arguments \cite{Clozel-Inventiones}; in the CM case we use Mok's classification \cite{Mok}  of the discrete spectrum for unitary groups; 
and for the parallel weight case it is an easy generalization of the constructive proof in Clozel~\cite{Clozel-Duke}. Finally, we add some remarks on an endoscopic stratification of inner cohomology suggested by results on arithmeticity as in Gan--Raghuram~\cite{Gan-Raghuram}


\bigskip
\section{Statements of the main results}

\medskip
\subsection{Cuspidal cohomology of $\GL_n$} 
\label{sec:cuspidal-cohomology}
We need a somewhat lengthy notational preparation; the reader familiar with such objects can directly jump to 
(\ref{eqn:cusp-coh-defn-2}). 
Let $F$ be a number field of degree $d=[F:\Q]$ with ring of integers $\O$. For any place $v$ we write $F_v$ for the topological completion of $F$ at $v$. Let $S_{\infty}(F) = S_\infty$ be the set of archimedean places of $F$. Let $S_\infty = S_r \cup S_c$, where $S_r$ (resp., $S_c$) is the set of real (resp., complex) places. Let $\cE_F = {\rm Hom}(F,\C)$ be the set of all embeddings of $F$ as a field into $\C$. There is a canonical surjective map $\cE_F \to S_\infty$, which is a bijection on the real embeddings and real places, and identifies a pair of complex conjugate embeddings $\{\iota_v, \bar{\iota}_v\}$ with the  complex place $v$. For each $v \in S_r$, we fix an isomorphism 
$F_v \cong \R$ which is canonical. Similarly for $v \in S_c$, we fix $F_v \cong \C$ given by (say) $\iota_v$; this choice is not canonical. 
Let $r_1 = |S_r|$ and $r_2 = |S_c|$; hence $d = r_1 + 2r_2.$ 
If $v \notin S_{\infty}$, we let $\O_v$ be the ring of integers of $F_v$, and $\wp_v$ it's  unique maximal ideal. 
Moreover, $\A_F$ denotes the ring of ad\`eles of $F$ and $\A_{F,f}$ its finite part. The group of id\`eles of $F$ will be denoted $\A_F^\times$ and similarly, $\A_{F,f}^\times$ is the group of finite id\`eles. We will drop the subscript $F$ when talking about 
$\Q$. Hence, $\A$ is $\A_\Q$, etc. 

\smallskip

The algebraic group ${\rm GL}_n/F$ will be denoted as $\underline{G}_n$, and we put 
$G_n = {\rm Res}_{F/\Q}(\underline{G}_n)$; an $F$-group will be denoted by 
an underline and the corresponding $\Q$-group via Weil restriction of scalars will be denoted without the underline; hence for any $\Q$-algebra $A,$ the 
group of $A$-points of $G_n$ is $G_n(A) = \underline{G}_n(A \otimes_\Q F)$. 
Let $\underline{B}_n = \underline{T}_n \underline{U}_n$ stand for the standard Borel subgroup of $\underline{G}_n$ of all upper triangular matrices, where $\underline{U}_n$ is the unipotent radical of $\underline{B}_n$, and $\underline{T}_n$ the diagonal torus. The center of $\underline{G}_n$ will be denoted by $\underline{Z}_n$. These groups define the corresponding $\Q$-groups $G_n \supset B_n = T_nU_n \supset Z_n$. Observe that $Z_n$ is not $\Q$-split, and we let $S_n$ be the maximal $\Q$-split torus in $Z_n$; we have 
$S_n \cong {\mathbb G}_m$ over $\Q.$

Note that 
$$
G_{n,\infty} := G_n(\R) = \underline{G}_n(F \otimes_\Q \R) = \prod_{v \in S_{\infty}} \GL_n(F_v) \cong \prod_{v \in S_r} \GL_n(\R) \times \prod_{v \in S_c} \GL_n(\C).
$$ 
We have $Z_n(\R) \simeq  \prod_{v \in S_r} \R^\times \times \prod_{v \in S_c} \C^\times.$ The subgroup $S_n(\R)$ is 
$\R^\times$ diagonally embedded in $Z_n(\R).$
Let $C_{n,\infty} =  \prod_{v \in S_r} \Orth(n) \times \prod_{v \in S_c} \Unit(n)$ be the maximal compact subgroup of $G_n(\R)$, and let 
$K_{n,\infty} \ =  \ Z_n({\mathbb R}) C_{n,\infty} \ = \ Z_n({\mathbb R})^0C_{n,\infty}.$ 
Let $K_{n,\infty}^0$ be the topological connected component of $K_{n,\infty}$. 
For a real Lie group $G$, we denote its Lie algebra by $\g^0$ and the complexified Lie algebra by $\g,$ i.e., 
$\g = \g^0 \otimes_\R \C.$ If $G = \GL_n(\R)$ then $\g^0 = \gl_n(\R)$ and $\g = \gl_n(\C)$, 
on the other hand, if $G$ stands for the real Lie group $\GL_n(\C)$ then $\g^0 = \gl_n(\C)$ as a Lie algebra over $\R$, 
and $\g = \gl_n(\C) \otimes_\R \C.$  With this notational scheme, we have $\g_n$, 
$\mathfrak{b}_n$, $\mathfrak{t}_n$ and $\k_n$ denoting the complexified Lie algebras of $G_n(\R)$, 
$B_n(\R)$, $T_n(\R)$ and $K_{n,\infty}^0$ respectively. For example, 
$\g_n = \prod_{v \in S_r} \gl_n(\C)  \times \prod_{v \in S_c}( \gl_n(\C) \otimes_\R \C).$

\smallskip

Consider $T_n(\R) = \underline{T}_n(F \otimes {\mathbb R}) \cong \prod_{v \in S_\infty} \underline{T}_n(F_v)$. 
We let $X^*(T_n \times \C)$ stand for the group of all algebraic characters of $T_n \times \C$, and let 
$X^+(T_n \times \C)$ stand for all those characters in $X^*(T_n \times \C)$ 
which are dominant with respect to $B_n$. A weight $\mu \in 
X^+(T_n \times \C)$ may be described as follows: 
$\mu \ = \ (\mu^\iota)_{\iota \in \cE_F},$ also sometimes written as $\mu = (\mu^v)_{v \in S_{\infty}},$ where $\mu^v = (\mu^{\iota_v}, \mu^{\bar{\iota}_v})$ for $v \in S_c,$ and furthermore 
\begin{itemize}
\item 
for $v \in S_r$, we have $\mu^v = (\mu^v_1,\dots, \mu^v_n)$, $\mu^v_i \in {\mathbb Z}$, 
$\mu^v_1 \geq \cdots \geq \mu^v_n$, 
and the character $\mu^v$ sends $t = {\rm diag}(t_1,\dots,t_n) \in \underline{T}_n(F_v)$ to 
$\prod_i t_i^{\mu^v_i},$  and  
\smallskip
\item if $v \in S_c$ then $\mu^v$ is a pair $(\mu^{\iota_v}, \mu^{\bar{\iota}_v})$, with $\mu^{\iota_v} = (\mu^{\iota_v}_1,\dots, \mu^{\iota_v}_n)$, 
 $\mu^{\iota_v}_i \in {\mathbb Z}$, $\mu^{\iota_v}_1 \geq \cdots \geq \mu^{\iota_v}_n$; likewise
$\mu^{\bar{\iota}_v} = (\mu^{\bar{\iota}_v}_1,\dots, \mu^{\bar{\iota}_v}_n)$ and  
$\mu^{\bar{\iota}_v}_1 \geq \cdots \geq \mu^{\bar{\iota}_v}_n;$  
the character $\mu^v$ is given by sending $t = {\rm diag}(z_1,\dots,z_n) \in \underline{T}_n(F_v)$ to 
$\prod_{i=1}^n z_i^{\mu^{\iota_v}_i} \bar{z}_i^{\mu^{\bar{\iota}_v}_i},$ where $\bar{z}_i$ is the conjugate of $z_i$.  
\end{itemize}

\smallskip

For $\mu \in X^+(T_n \times \C)$, define a finite-dimensional complex 
representation $(\rho_{\mu}, \M_{\mu, \C})$ of $G_n(\R)$ as follows: 
For $v \in S_r$, let $(\rho_{\mu^v}, \M_{\mu^v, \C})$ 
be the irreducible complex representation of $G_n(F_v) = \GL_n({\mathbb R})$ with 
highest weight $\mu_v$.
For $v \in S_c$, let $(\rho_{\mu^v}, \M_{\mu^v, \C})$ be the complex representation of the real algebraic group 
$G(F_v) = \GL_n({\mathbb C})$ defined as 
$\rho_{\mu^v}(g) = \rho_{\mu^{\iota_v}}(g) \otimes \rho_{\mu^{\bar{\iota}_v}}(\overline{g});$ here
$\rho_{\mu^{\iota_v}}$ (resp., $\rho_{\mu^{\bar{\iota}_v}}$)
is the irreducible representation of the complex group $\GL_n(\C)$ with highest weight 
$\mu^{\iota_v}$ (resp., $\mu^{\bar{\iota}_v}$). 
Now we let $\rho_{\mu} = \otimes_{v \in S_{\infty}} \rho_{\mu^v}$ which acts on 
$\M_{\mu, \C} = \otimes_{v \in S_{\infty}} \M_{\mu^v, \C}.$

\smallskip

Let $K_f$ be an open compact subgroup of $G_n(\A_f) = \GL_n(\A_{F,f})$. Let us write $K_f = \prod_p K_p$ where each $K_p$ is an open compact subgroup of $G_n(\Q_p)$ and for almost all $p$ we have $K_p = \prod_{v | p} \GL_n(\O_v)$. 
Define the double-coset space 
$$
S^{G_n}_{K_f} \ = \ G_n(\Q) \backslash G_n(\A) /K_{n,\infty}^0 K_f 
\ = \ \GL_n(F) \backslash \GL_n(\A_F) /K_{n,\infty}^0 K_f.
$$
Given a dominant-integral weight $\mu \in X^+(T_n \times \C)$ we get a sheaf $\tM_{\mu, \C}$ of $\C$-vector spaces on $S^{G_n}_{K_f}.$ (See, for example, \cite[Sect.\,2.3.3]{Raghuram-Forum}.)
We are interested in the 
sheaf cohomology groups 
$H^{\bullet}(S^{G_n}_{K_f} , \tM_{\mu, \C}^{\sf v}).$
Here $\tM_{\mu, \C}^{\sf v}$ is the sheaf attached to the contragredient representation $\M_{\mu,\C}^{\sf v}$  of 
$\M_{\mu, \C}$. This dualizing is only for convenience. Let $\omega_{\rho_\mu}$ be the central character of $\rho_\mu.$ 
For the sheaf $\tM_{\mu, \C}^{\sf v}$ or $\tM_{\mu, \C}$ to be nonzero we need the condition: 
\begin{equation}
\label{eqn:sheaf-nonzero}
\mbox{``The central character $\omega_{\rho_\mu}$ is trivial on $Z_n(\Q) \cap K_{n,\infty}^\circ K_f$."}
\end{equation}
Henceforth, we will assume this condition on $\mu$ and $K_f.$ (This can be a nontrivial condition even in simple situations: take, for example, $F=\Q$, $n=2$, $K_f = \prod_p \GL_2(\Z_p)$, and $\mu = (a,b)$ with integers $a \geq b$; then this condition boils down to $(-1)^{a+b} = 1;$ but by taking $K_f$ slightly deep enough we can ensure 
$Z_2(\Q) \cap K_{2,\infty}^\circ K_f$ is trivial and so the condition vacuously holds.) For more details, the reader is referred to Harder~\cite[(1.1.3)]{Harder-Inventiones}. 
Passing to the limit over all open compact subgroups $K_f$ and let 
$H^{\bullet}(S^{G_n}, \tM_{\mu,\C}^{\sf v}) := \varinjlim_{K_f} H^{\bullet}(S^{G_n}_{K_f} , \tM_{\mu,\C}^{\sf v}).$
There is an action of $\pi_0(G_{n,\infty}) \times G_n(\A_f)$ on $H^{\bullet}(S^{G_n}, \tM_{\mu,E}^{\sf v})$, and 
the cohomology of $S^{G_n}_{K_f}$ is obtained by taking invariants under $K_f$, i.e., 
$H^{\bullet}(S^{G_n}_{K_f} , \tM_{\mu,E}^{\sf v}) = H^{\bullet}(S^{G_n}, \tM_{\mu,E}^{\sf v})^{K_f}.$ 
We can compute the above sheaf cohomology via the de~Rham complex, and then reinterpreting the de~Rham complex in terms of the complex computing relative Lie algebra cohomology, we get the isomorphism: 
$$
H^{\bullet}(S^{G_n}, \tM_{\mu,\C}^{\sf v})  \ \simeq \ 
H^{\bullet}(\g_n, K_{n,\infty}^0; \  C^{\infty}(G_n(\Q)\backslash G_n(\A)) \otimes \M_{\mu,\C}^{\sf v}) .
$$

\smallskip

The inclusion $C^{\infty}_{\rm cusp} (G_n(\Q)\backslash G_n(\A)) \hookrightarrow C^{\infty}(G_n(\Q)\backslash G_n(\A))$ of the space of smooth cusp forms  in the space of all smooth functions induces, via results of Borel \cite{borel-duke},  an injection in cohomology; this defines cuspidal cohomology: 
\begin{equation}
\label{eqn:cusp-coh-defn-1}
\xymatrix{
H^{\bullet}(S^{G_n}, \tM_{\mu,\C}^{\sf v}) 
\ar[rr] & &
H^{\bullet}(\g_n, K_{n,\infty}^0; C^{\infty}(G_n(\Q)\backslash G_n(\A)) \otimes \M_{\mu,\C}^{\sf v})  \\
H^{\bullet}_{\rm cusp}(S^{G_n}, \tM_{\mu,\C}^{\sf v}) \ar@{^{(}->}[u]
\ar[rr] 
& & 
H^{\bullet}(\g_n, K_{n, \infty}^0; C^{\infty}_{\rm cusp}(G_n(\Q)\backslash G_n(\A)) \otimes \M_{\mu,\C}^{\sf v}) 
\ar@{^{(}->}[u]
}
\end{equation}
With level structure $K_f$ we have:  
\begin{equation}
\label{eqn:cusp-coh-defn-2}
H^{\bullet}_{\rm cusp}(S^{G_n}_{K_f}, \tM_{\mu,\C}^{\sf v})   \simeq  
H^{\bullet}(\g_n, K_{n,\infty}^0; \  C^{\infty}_{\rm cusp}(G_n(\Q)\backslash G_n(\A))^{K_f} \otimes \M_{\mu,\C}^{\sf v}) .
\end{equation}

\smallskip

A fundamental problem concerning cuspidal cohomology for $\GL_n/F$ is: 
\begin{prob}
\label{props:strong}
Classify all $(F,n,\mu,K_f)$, subject to (\ref{eqn:sheaf-nonzero}), 
for which $H^{\bullet}_{\rm cusp}(S^{G_n}_{K_f}, \tM_{\mu,\C}^{\sf v}) \neq 0.$ 
\end{prob}
This problem includes classical situations such as: take $F = \Q$ and $n=2$, and given integers $N \geq 1$ and $k \geq 2$ is there a holomorphic cusp form of weight $k$ for $\Gamma_1(N) \subset \SL_2(\Z)$? One may relax the dependence on an explicit level structure $K_f$ and ask for a solution of 
the weaker 
\begin{prob}
\label{prob:weak}
Classify all $(F,n,\mu)$ for which $H^{\bullet}_{\rm cusp}(S^{G_n}_{K_f}, \tM_{\mu,\C}^{\sf v}) \neq 0$ for some $K_f.$
\end{prob}

\medskip

{\it The purpose of this article is to provide a solution of the weaker problem for $F$ a totally real field, or a totally imaginary quadratic extension of a totally real field, or for a general number field $F$ and for a parallel weight $\mu$; see below for a definition of a parallel weight.}

\subsection{Necessary conditions for nonvanishing of cuspidal cohomology}

Using the usual decomposition of the space of cusp forms into a direct sum of cuspidal automorphic representations, we get the following fundamental decomposition of 
$\pi_0(G_n(\R)) \times G_n(\A_f)$-modules: 
\begin{equation}
\label{eqn:cusp-coh}
H^{\bullet}_{\rm cusp}(S^{G_n}, \tM_{\mu,\C}^{\sf v}) \ = \ 
\bigoplus_{\Pi} H^{\bullet}(\g_n, K_{n,\infty}^0;  \Pi_{\infty} \otimes  \M_{\mu,\C}^{\sf v}) \otimes \Pi_f.
\end{equation}
We say that {\it $\Pi$ contributes to the cuspidal cohomology of $G_n$ with coefficients in $\M_{\mu,\C}^{\sf v}$}, and 
we write $\Pi \in {\rm Coh}(G_n, \mu^{\sf v})$,  if $\Pi$ has a nonzero contribution to the above decomposition. Equivalently, if $\Pi$ is a cuspidal automorphic representation whose representation at infinity $\Pi_{\infty}$ after twisting by $\M_{\mu,\C}^{\sf v}$ has nontrivial relative Lie algebra cohomology. With a level structure $K_f$, (\ref{eqn:cusp-coh}) takes the form: 
\begin{equation}
\label{eqn:cusp-coh-level}
H^{\bullet}_{\rm cusp}(S^{G_n}_{K_f}, \tM_{\mu,\C}^{\sf v}) \ = \ 
\bigoplus_{\Pi} H^{\bullet}(\g_n, K_{n,\infty}^0;  \Pi_{\infty} \otimes  \M_{\mu,\C}^{\sf v}) \otimes \Pi_f^{K_f}. 
\end{equation}
We write $\Pi \in {\rm Coh}(G_n, \mu^{\sf v}, K_f)$ if a cuspidal automorphic representation $\Pi$ contributes nontrivially to (\ref{eqn:cusp-coh-level}). Note that each ${\rm Coh}(G_n, \mu^{\sf v}, K_f)$ is a finite set, and clearly, 
${\rm Coh}(G_n, \mu^{\sf v}) = \cup_{K_f} {\rm Coh}(G_n, \mu^{\sf v}, K_f).$ 

\smallskip

Suppose $\Pi \in \Coh(G_n, \mu^{\sf v}).$ The fact that $\mu^{\sf v}$ supports cuspidal cohomology places some restrictions on $\mu.$ First of all, essential unitarity of $\Pi$, and in particular of $\Pi_\infty$ gives, via Wigner's Lemma, essential self-duality of $\mu$: there is an integer ${\sf w}(\mu)$ such that 
\begin{enumerate}
\item For $v \in S_r$ and $1 \leq i \leq n$ we have $\mu^{\iota_v}_i + \mu^{\iota_v}_{n-i+1} = {\sf w}(\mu);$ 
\smallskip
\item For $v \in S_c$ and $1 \leq i \leq n$ we have $\mu^{\bar{\iota}_v}_i + \mu^{\iota_v}_{n-i+1} = {\sf w}(\mu).$ 
\end{enumerate}
We will call such a weight $\mu$ as a {\it pure weight} and call ${\sf w}(\mu)$ the {\it purity weight of $\mu.$} Let 
$X^+_0(T_n)$ denote the set of dominant integral pure weights. Applying 
Clozel \cite[Theorem 3.13]{Clozel-Motives}, we get 
$$
\Pi \in \Coh(G_n, \mu^{\sf v}) \ \Longrightarrow \  {}^\sigma\Pi \in \Coh(G_n, {}^\sigma\!\mu^{\sf v}), \quad \forall 
\sigma \in \autc,
$$
where, ${}^\sigma\!\mu \in X^*(T_n \times \C)$ is defined as: 
${}^\sigma\!\mu \ = \ ({}^\sigma\!\mu^\iota)_{\iota \in \cE_F}$ with  
${}^\sigma\!\mu^\iota := \mu^{\sigma^{-1} \circ \iota}.$ The reader is referred to \cite{Clozel-Motives} or \cite{Raghuram-Forum} for a definition of ${}^\sigma\Pi.$
In particular, ${}^\sigma\!\mu$ also satisfies the purity conditions (1) and (2) above. Note that 
${\sf w}(\mu) = {\sf w}({}^\sigma\!\mu).$

\begin{defn}\label{def:strongly-pure}
Let $\mu \in X^+_0(T_n)$ be a pure dominant integral weight. 
We say $\mu$ is {\it strongly pure} if 
${}^\sigma \!\mu$ is pure for all $\sigma \in \autc.$ Let $X_{00}^+(T_n)$ stand for the set of dominant integral strongly pure weights. 
\end{defn}

For any $F$, we have the following inclusions $X^+_{00}(T_n) \subset X^+_0(T_n) \subset X^+(T_n) \subset X^*(T_n)$ and in general they are all strict inclusions. 
If $F$ is a totally real field or a CM field (totally imaginary quadratic extension of a totally real field) 
then $\mu$ is pure if and only if $\mu$ is strongly pure. For any number field, one may see that there are strongly pure weights: take an integer $b$ and integers 
$a_1 \geq a_2 \geq \cdots \geq a_n$ such that $a_j + a_{n-j+1} = b;$ now for 
each $\iota \in \cE_F$ put $\mu^\iota = (a_1,\dots, a_n)$;  then $\mu$ is strongly pure with ${\sf w}(\mu) = b;$
such a weight may be called a {\it parallel weight}. We formulate the following conjecture as a possible answer to 
Problem~\ref{prob:weak}. 

\smallskip

\begin{con}
Let $\mu \in X^+(T_n \times \C)$ be a dominant integral weight. Then 
$$H^{\bullet}_{\rm cusp}(S^{G_n}_{K_f}, \tM_{\mu,\C}^{\sf v}) \neq 0 \ \mbox{for some $K_f$} \ \ \iff \ \  
\mu \in X^+_{00}(T_n \times \C). 
$$
\end{con}

\smallskip

We give a brief glimpse, without pretending to be exhaustive, into available results in the literature. 
Borel, Labesse and Schwermer \cite[Theorem 11.3]{Borel-Labesse-Schwermer} proved nonvanishing of cuspidal cohomology for $G = \SL_n$ over any number field which is an extension via cyclic prime degree extensions of a totally real number field, and for the trivial coefficient system, i.e., when $\mu = 0;$ their proof involved studying Lefschetz numbers for rational automorphisms of  $G.$ Barbasch and Speh \cite[\S XI]{Barbasch-Speh} considered $\GL_n/\Q$ and a coefficient system  with certain technical restrictions on the weight $\mu$ (which in fact exclude many pure weights $\mu$) and proved nonvanishing of cuspidal cohomology via an application of trace formula for certain Lefschetz functions.  Clozel \cite[Theorem 4]{Clozel-Duke} gave a constructive proof of the nonvanishing of cuspidal cohomology for $\GL_{2n}$ over any number field for the trivial coefficient system using automorphic induction.

\medskip
\subsection{The main results of this article on cuspidal cohomology of $\GL(n)$} 

\medskip

\begin{theorem}
\label{thm:main} 
Take an integer $N \geq 2.$ 
Let $F_0$ be a totally real field extension of $\Q$, and we take an extension $F/F_0$ to be either
\begin{enumerate}
\item $F = F_0$, i.e., $F$ itself is a totally real field; or 
\item $F$ is a totally imaginary quadratic extension over $F_0.$ 
\end{enumerate}
Let $G = G_N = {\rm Res}_{F/\Q}(\GL_N/F).$ Let $\mu \in X^+_0(T_N \times \C)$ be a pure dominant integral weight with purity ${\sf w}(\mu) = 0.$ In case (2), assume furthermore that $\mu$ is trivial on the center $Z_N$, i.e., for all $\iota \in \cE_F$ suppose that $\mu^\iota_1 + \dots + \mu^\iota_N = 0.$ Then
$$
H^\bullet_{\rm cusp}(S^G, \M_{\mu, \C}^{\sf v}) \ \neq \ 0.
$$
\end{theorem}

The proof involves endoscopic transfer from certain classical groups and breaks up into the following sub-cases: 
\begin{enumerate}
\item[(1a)] $F$ is totally real, and $N = 2n+1$ is odd; 
\item[(1b)] $F$ is totally real, and $N = 2n$ is even;  
\item[(2)] $F$ is a totally imaginary quadratic extension of a totally real $F_0.$ 
\end{enumerate}

In case (1a), we transfer from the group $G' = \Sp(2n)/F.$ To begin, we transfer the weight $\mu$ to a weight $\mu'$ on $G'.$  Then, using results on limit multiplicities due to Clozel~\cite{Clozel-Inventiones}, 
we produce a cuspidal representation $\pi'$ of $G'$ with a discrete series representation at infinity which is cohomological with respect to $\mu'$. Furthermore, one can arrange for $\pi'$ to have the
Steinberg representation at some finite place. Arthur's results ensure then that $\pi'$ corresponds to an Arthur 
parameter $\pi$ on $G$ which is cuspidal; the cuspidality of $\pi$ uses an observation of Magaard and Savin \cite{Magaard-Savin},  
and then one checks that indeed $\pi$ contributes to the above cohomology group.  
For the proof in 
case (1b) (resp., case (2)) we transfer from the split group $G' = \SO(2n+1)$ (resp., $G'$ a unitary group in $N$ variables). In case (1a) it might be possible to use the results of Weselmann's papers \cite{Weselmann-1} and 
\cite{Weselmann-2}.

\medskip

We may draw several inferences from the above proof. As summarized in \cite[Section 3.1]{Grobner-Raghuram}, a cuspidal representation $\pi$ of $\GL_{2n}/F$ is a transfer from $\SO(2n+1)$ if and only if a partial exterior-square $L$-function $L^S(s, \wedge^2, \pi)$ has a pole at $s=1$ and this is so if and only if $\pi$ admits a Shalika model, which gives us

\begin{cor}
Let $F$ be a totally real field, and take $G = \GL(2n)/F.$ Let $\mu$ be a pure weight with purity $0.$ Then there exists a cuspidal automorphic representation $\pi$ of $G$ such that 
\begin{enumerate}
\item $\pi \in \Coh(G, \mu),$ i.e., it is cohomological with respect to $\mu$, and 
\item $\pi$ has a Shalika model, or equivalently, a partial exterior-square $L$-function $L^S(s, \wedge^2, \pi)$ has a pole at $s=1.$
\end{enumerate}
\end{cor}

A similar characterization for a $\pi$ on $\GL_{2n+1}/F$ being a transfer from $\Sp(2n)$ if and only if a partial 
symmetric-square $L$-function $L^S(s, {\rm Sym}^2, \pi)$ having a pole at $s=1$ gives us 

\begin{cor}
Let $F$ be a totally real field, and take $G = \GL(2n+1)/F.$ Let $\mu$ be a pure weight with purity $0.$ Then there exists a cuspidal automorphic representation $\pi$ of $G$ such that 
\begin{enumerate}
\item $\pi \in \Coh(G, \mu),$ i.e., it is cohomological with respect to $\mu$, and 
\item a partial symmetric-square $L$-function $L^S(s, {\rm Sym}^2, \pi)$ has a pole at $s=1.$ 
\end{enumerate}
\end{cor}

\medskip

Using an idea in Labesse and Schwermer~\cite{Labesse-Schwermer}, that at a finite place the Steinberg representation retains the property of being Steinberg upon base change, and at archimedean places the property of being cohomological is preserved under base-change, we get the following

\begin{cor}
Let $F$ be a totally real field and suppose that $\tilde F/F$ is a finite extension that is filtered by cyclic extensions of prime degrees. Let $\mu$ be a pure weight for $G = \GL(N)/F$ with purity $0.$ Define a weight $\tilde \mu$ for
$\tilde G = \GL(N)/\tilde F$ as
for any $\tilde \iota : \tilde F \to \C$, we let $\tilde \mu^{\tilde \iota} = \mu^\iota$ where $\iota = \tilde\iota|_F.$ (If $\mu$ is a parallel weight then so is $\tilde\mu.$) Then
$$
H^\bullet_{\rm cusp}(S^{\tilde G}, \M_{\tilde\mu, \C}^{\sf v}) \ \neq \ 0.
$$
\end{cor}
The conditions on $\tilde F/F$ are dictated by the main theorem on base change for $\GL_N$ due to Arthur and Clozel
\cite{Arthur-Clozel}.

\medskip

Next, we consider a general number field $F$. The following theorem is a generalization of the main theorem of 
Clozel~\cite{Clozel-Duke} which is proved by constructing cohomological cuspidal representations via automorphic induction. See also Ramakrishnan--Wang \cite[Appendix]{Ramakrishnan-Wang}.

\begin{theorem}
\label{thm:auto-induction}
Let $F$ be any number field, and take $G = \GL(2n)/F.$ Let $\mu$ be a parallel weight with purity ${\sf w}(\mu) = 0.$ 
Then
$$
H^\bullet_{\rm cusp}(S^G, \M_{\mu, \C}^{\sf v}) \ \neq \ 0.
$$
\end{theorem}

\bigskip
\section{Proof of Theorem~\ref{thm:main}}

\medskip
\subsection{Archimedean preliminaries}

\medskip
\noindent{\bf Case (1a) ($F$ is totally real, and $N = 2n+1$ is odd.)}
Fix a real place $v$ of $F$, and since $\mu$ is a pure weight with purity $0$, we have 
$$
\mu^v \ = \ 
\mu^v_1 \geq \mu^v_2 \geq \cdots \geq \mu^v_n \geq 0 \geq -\mu^v_n \geq \cdots
\geq -\mu^v_2 \geq -\mu^v_1.
$$ 
The place $v$ of $F$ is fixed, and for brevity, we will drop $v$ from the notation for $\mu^v$.

\smallskip

Consider the endoscopy group $G' = \Sp(2n)/F$ defined so that upper-triangular subgroup $B'$ in $G'$ is a Borel subgroup. 
The connected component of the $L$-group of $G'$ is ${}^L\!G'^\circ = \SO(2n+1,\C)$. The maximal compact subgroup $K'$ of $G'(\R) = \Sp(2n,\R)$ is isomorphic to $\Unit(n).$ Define a dominant integral weight $\mu'$ for $G'$ which at the place $v$ (dropped from the notation) is given by: 
$$
\mu':=  (\mu_1, \mu_2, \dots, \mu_n) \ = \ \sum_{i=1}^n \mu_i e_i, 
$$
where $e_i$ gives the $i$-th coordinate of a diagonal matrix. 
Let $\rho'$ be the half-sum of positive roots for $G'$, which is written as  
$\rho' = \sum \limits_{j=1}^{n} (n+1-j) e_j = (n, n - 1, \dots, 1).$ Let 
$$
\Lambda' \ = \ \mu' + \rho'  \ = \ (\mu_1+n, \ \mu_2+n-1, \dots, \ \mu_{n-1}+1, \ \mu_n).
$$ 
Thus $\Lambda'$ is a regular weight and using Harish-Chandra's classification theorem of discrete series representations (see, for example, Knapp~\cite{Knapp-book}), there exists a
discrete series representation $\pi' = \pi_{\Lambda'}$ of $G'(\R)$ whose infinitesimal character is $\chi_{\Lambda'}$. 
Let $\mathcal{M}_{\mu',\C}$ be the algebraic irreducible representation of $G'(\C)$ with highest weight $\mu'$.
From the well-known results on the cohomology of discrete series representations (see, for example, 
Borel--Wallach~\cite[Theorem II.5.3]{Borel-Wallach}), one knows 
that $\pi'$ is cohomological with respect to the coefficient system $\mathcal{M}_{\mu',\C}$ of $G',$ 
i.e., the relative Lie algebra cohomology $H^{\mathbf{\bullet}}(\mathfrak{g'}_{\infty}, K', \pi' \otimes \mathcal{M}^{\sf v}_{\mu'} )$ is non-zero; in fact it is nonzero only in the middle degree $\bullet = \tfrac12 {\rm dim}(G'(\R))/{\rm dim}(K').$

The shape of the Langlands parameter, denoted $\tau_{\Lambda'},$ of the discrete series representation 
$\pi_{\Lambda'}$ of $G'(\R)$ is well-known; we can deduce the following from \cite[Example 10.5]{Borel-corvallis}: 
$$
\tau_{\Lambda'} \ = \ 
\Ind_{\C^\times}^{W_\R}(\chi_{\ell_1}) \oplus \Ind_{\C^\times}^{W_\R}(\chi_{\ell_2}) \oplus \cdots \oplus 
\Ind_{\C^\times}^{W_\R}(\chi_{\ell_n}) \oplus {\rm sgn}^{\epsilon}, 
$$
where $\ell_1, \dots, \ell_n$ are positive integers and the first $n$-summands are irreducible $2$-dimensional representations of the Weil group $W_\R$ of $\R$ induced from characters of $\C^\times$ with the character $\chi_{\ell_j}$ sending $z = re^{i \theta} \in \C^\times$  to 
$r^{i \ell_j \theta} = (z/\bar{z})^{\ell_j/2}$ with the proviso that each summand be of orthogonal type forcing 
each $\ell_j$ to be even, and the determinant of the parameter is $1$ so as to have image inside ${}^L\!G'^\circ$ 
forcing the last summand to be ${\rm sgn}^n$. Furthermore, the relation between the integers $\ell_j$ and the 
weight $\Lambda'$ is captured by 
$$
\tau_{\Lambda'}|_{\C^\times} \ = \ z^{\Lambda'} \bar{z}^{-\Lambda'}. 
$$
(Here we have tacitly used that $\Lambda'$, which is a character of a maximal torus $T'$ of $G',$ is also, by definition, a co-character of the dual ${}^L\!T'^\circ \subset {}^L\!G'^\circ$ justifying the above notation.) 
If we put $\ell = (\ell_1,\dots,\ell_n)$ then we get $\ell  \ = \ 2\Lambda',$ i.e., 
$$
(\ell_1,\dots,\ell_n) \ = \ (2\mu_1+2n, \ 2\mu_2+2n-2, \dots, \ 2\mu_{n-1}+2, \ 2\mu_n). 
$$

Let $\pi_\mu$ denote the Langlands transfer of $\pi'$ to an irreducible representation of $\GL_{2n+1}(\R)$, where the transfer is mitigated by the Langlands parameter of $\pi'$ being that of $\pi_\mu$ via the standard embedding 
${}^L\!G'^\circ = \SO(2n+1,\C) \subset \GL(2n+1,\C) = {}^L\!G^\circ$. Using the local Langlands correspondence for 
$\GL_N(\R)$ (see, for example, Knapp~\cite{Knapp}), we can deduce 
$$
\pi_\mu  \ = \ \Ind^{G}_{P(2,2,\dots, 2,1)} \left( D_{ \ell_1} \otimes
 D_{ \ell_2} \otimes \cdots \otimes  D_{ \ell_n} \otimes {\rm sgn}^{n} \right), 
 $$
 where, for any integer $l$, we denote by $D_l$ the discrete series representation of $\GL_2(\R)$ as normalized in 
 \cite[3.1.3]{Raghuram-Tanabe}. It is well-known (\cite[Lemme 3.14]{Clozel-Motives}) that 
 $$
 H^\bullet(\gl_N, \R^\times_+\SO(N); \, \pi_\mu \otimes \M_{\mu,\C}) \ \neq \ 0.
 $$

\medskip
\noindent{\bf Case (1b) ($F$ is totally real, and $N = 2n$ is even.)} 
In this case we take $G' = \SO(2n+1)/F$ which is the split orthogonal group in $2n+1$ variables. 
We have $G'(\R) = \SO(n,n+1).$ The maximal compact subgroup $K'$ of $G'(\R)$ is isomorphic to 
${\rm S}({\rm O}(n) \times {\rm O}(n+1))$. The connected component of the $L$-group is 
${}^L\!G'^\circ = \Sp(2n,\C).$ Fix a real place $v$ of $F$ and we drop it from the notations. We have 
$$
\mu = (\mu_1 \geq \mu_2 \geq \cdots \geq \mu_n \geq -\mu_n \geq \cdots
\geq -\mu_2 \geq -\mu_1)
$$
be the given dominant integral weight with purity weight $0$. Following the argument as in Case (1) above, we put
$$
\begin{array}{lll}
\mu' & = & (\mu_1, \mu_2, \dots, \mu_n), \quad {\rm and} \\ 
& & \\ 
\Lambda' & = & \mu' + \rho' \ = \ (\mu_1+n -\frac{1}{2}, \ \mu_2+n - \frac{3}{2}, \ \cdots, \ \mu_{n-1}+ \frac{3}{2}, \ 
\mu_n+\frac{1}{2}). 
\end{array}
$$
Consider the discrete series representation $\pi' = \pi_{\Lambda'}$ with infinitesimal character given by $\Lambda'.$ The Langlands parameter $\tau_{\Lambda'}$ of $\pi_{\Lambda'}$ has the form  
$$
\tau_{\Lambda'} \ = \ \Ind_{\C^\times}^{W_\R}(\chi_{\ell_1}) \oplus \Ind_{\C^\times}^{W_\R}(\chi_{\ell_2}) \oplus \cdots \oplus 
\Ind_{\C^\times}^{W_\R}(\chi_{\ell_n}),
$$ 
with all the $\ell_j$ being odd positive integers. The infinitesimal character of the discrete series is seen in terms of the exponents of the parameter restricted to $\C^\times$ giving us $\ell  \ = \ 2\Lambda',$ or that 
$$
(\ell_1,\dots,\ell_n) \ = \ (2\mu_1+2n - 1, \ 2\mu_2 + 2n - 3, \ \cdots, 2\mu_{n-1}+3, \ 2\mu_n+1).
$$
Via the local Langlands correspondence for $\GL_{2n}(\R)$ we get that $\pi'$ transfers to $\pi_\mu$ given by
$$
\pi_\mu \ = \ 
\Ind^{G}_{P(2,2,\dots, 2)} \left( D_{ \ell_1} \otimes
 D_{ \ell_2} \otimes \cdots \otimes  D_{ \ell_n} \right), 
 $$
which has the property that $H^\bullet(\gl_N, \R^\times_+\SO(N); \, \pi \otimes \M_{\mu,\C}) \neq 0.$

\medskip
\noindent{\bf Case (2) ($F_0$ is totally real and $F$ a totally imaginary
quadratic extension of $F_0$.)}
Let $\mu$ be a pure dominant integral weight for $G$ with purity weight $0.$ For a (complex) place $v$ of $F$, we have 
$\mu^v = (\mu^{\iota_v}, \mu^{\bar{\iota}_v})$. The place $v$ is fixed and we drop it from the notations. We have 
$\mu^{\iota} = (\mu_1 \geq \mu_2 \geq \cdots \geq \mu_{n})$ and 
$\mu^{\bar{\iota}} = (\mu_1^{\ast} \geq \mu_2^{\ast}\geq \cdots \geq \mu_{n}^{\ast})$ where purity, with the purity-weight being $0,$ implies $\mu_j^{\ast} = - \mu_{n-j+1}$.
In this case, we let $G'$ be the quasi-split unitary group over $F$; more precisely, define the matrix
$$
[\Phi] \ = \ \Phi_{i,j}, \quad {\rm where} \quad \Phi_{i,j} = (-1)^{i+1} \delta(j, n-i+1). 
$$
This gives a Hermitian form $\Phi$ on an $n$-dimensional $F$-vector space (with respect to a fixed basis) via
$\Phi(x,y) = {}^t\!x \cdot [\Phi] \cdot \bar{y},$ for $x = (x_1,\dots, x_n), y = (y_1,\dots, y_n) \in F^n$ and 
$\bar{y} = (\bar{y_1},\dots,\bar{y_n})$ being induced by the nontrivial Galois automorphism of $F/F_0.$ Denote 
$\Unit(\Phi)$ the corresponding unitary group over $F_0$, whose $F_0$-points consists of all $g \in \GL_n(F)$ such that 
${}^t\!g \cdot [\Phi] \cdot \bar{g} = [\Phi]$. Let $G'/F_0 := \Unit(\Phi).$ 
Note that $G'(F_0 \otimes \R)$ is a product of copies of $\Unit(\frac{n}{2}, \frac{n}{2})$ if $n$ is even, and is a product of  
copies of $\Unit( \frac{n+1}{2}, \frac{n-1}{2})$ if $n$ is odd. At $v$ (dropped from the notation) we put 
$$\mu' \ = \ (\mu_1, \mu_2, \dots, \mu_n).$$
The half sum $\rho'$ of positive roots is given by 
\[ 
\rho' = \frac{1}{2} \sum \limits_{1 \leq i < j \leq n} (e_i - e_j)=
 \left( \frac{n-1}{2},~ \frac{n-3}{2}, ~\cdots,  \frac{n-2i+1}{2},~ \cdots \frac{1-n}{2} \right).
 \]
 Once again define 
  \[ 
  \Lambda' \ : = \ \mu' + \rho' \ = \ 
  \left( \mu_1 + \frac{n-1}{2},~ \mu_2 + \frac{n-3}{2}, ~\dots, \ \mu_n + \frac{1-n}{2} \right), 
  \]
  which parametrizes a discrete series representation $\pi' = \pi_{\Lambda'}$ of 
  $\Unit(\lceil \frac{n}{2} \rceil, \lfloor \frac{n}{2} \rfloor).$ One may see, either from 
  \cite[Section 2]{Gan-Gross-Prasad} or from \cite[Example 10.5]{Borel-corvallis},  
 that the Langlands parameter $\tau_{\Lambda'}$ of $\pi_{\Lambda'}$ has the form  
  \[ 
  \tau_{\Lambda'} \ = \ z^{a_1} \bar{z}^{-a_1} \ \oplus \ z^{a_2} \bar{z}^{-a_2} \ \oplus \cdots \oplus \ z^{a_n} \bar{z}^{-a_n}, 
   \]
  where, once again, comparing the infinitesimal character with the exponents of the parameter we get  
  $(a_1,\dots,a_n) = 2 \Lambda'$, i.e., 
  $$
  (a_1,\dots, a_n) \ = \  \left( 2\mu_1 + n-1, \ 2\mu_2 + n-3,  \ \dots, \ 2\mu_n + 1-n\right). 
  $$
  
  The transfer $\pi_\mu$ of $\pi'$ to a representation $\GL(n,\C)$ is given by: 
  \[ 
  \pi_\mu \ =  \ \Ind^{GL(n,\C)}_{B(\C)}(z^{a_1} \bar{z}^{-a_1} \otimes \cdots \otimes \ z^{a_n} \bar{z}^{-a_n}), 
  \]
which has the property that $H^\bullet(\gl_n, \C^\times\Unit(n); \, \pi_\mu \otimes \M_{\mu,\C}) \neq 0.$  
 
\medskip

The above discussion in all the three cases, together with K\"unneth theorem for relative Lie algebra cohomology, gives the following

\begin{prop}
\label{prop:cohomology}
Let $F$ and $G$ be as in Theorem~\ref{thm:main}. Let $\mu$ be a pure dominant integral weight with purity $0.$ 
Define $G'$ as: 
\begin{enumerate}
\item[(1a)] if $F=F_0$ is totally real, and $N = 2n+1$ is odd; take $G' = \Sp(2n)/F;$ 
\item[(1b)] if $F=F_0$ is totally real, and $N = 2n$ is even;  take $G' = \SO(2n+1)/F$ (the split group); 
\item[(2)] if $F$ is a totally imaginary quadratic extension of a totally real $F_0;$ take $G' = \Unit(\Phi)/F_0$.   
\end{enumerate}
In each case, we define a dominant integral weight $\mu'$ for $G'$, and put $\Lambda' = \mu' + \rho'$. The discrete series representation 
$\pi_{\Lambda'}$ of $G'(F_0 \otimes \R)$ with infinitesimal character given by $\Lambda'$ has the property that it transfers to a representation $\pi_\mu$ of $G(\R) = \GL_N(F \otimes \R)$ that has nontrivial relative Lie algebra cohomology after twisting by $\M_{\mu,\C}.$ 
\end{prop}

\medskip
\subsection{Consequences of embedding theorems and global transfer}

\subsubsection{\bf Clozel's result on globalizing local discrete series representations}
In this section we describe a result of Clozel (\cite[Section 4.3]{Clozel-Inventiones})  on limit multiplicity of discrete series. Let $\cG$ be a semi-simple connected group defined over a number field $F.$  Let $v_0$ be a place such that
$\cG_{v_0}$ has supercuspidal representations. Let $S$ be a finite set of places containing the archimedean ones
and such that $v_0 \notin S$.
Let $S'$ be a finite set of finite places disjoint from $S$. Let $K_{S'}$
 be a compact open subgroup of $\cG_{S'} = \prod \limits_{v \in S'} \cG_p$.
Fix $K$, a compact open subgroup of $\prod \limits_{p \in S \cup S' \cup p_{0}} \cG_v$ and
 let $\mathcal{L}^{K_0 \times K_{S'} \times K}$ be the  $K_0 \times K_{S'} \times K$-invariant
functions in the space $L^{2}_{\text{cusp}}~( \cG(\Q) \backslash \cG(\mathbb{A}))$ of 
cusp forms on $\cG$. Let $\delta_S$ be a discrete series representation of the group $\cG_S$.

\begin{theorem}
\label{clozel limit multiplicity}
\[ 
\liminf \limits_{K_{S'} \rightarrow 1} v(K_{S'}) ~\text{mult}~(\delta_{S}, \mathcal{L}^{K_0 \times K_{S'} \times K}) \geq c
\]
where $c$ is a positive real number.
\end{theorem}

The above theorem is applicable to the simple groups $\SO(2n+1)$ or $\Sp(2n)$, however, it is not directly applicable to $\Unit(\Phi)$. It is for this reason, that in Theorem~\ref{thm:main}, in Case (2), we take $\mu$ to be trivial on the center, and so the discrete series representation $\pi_{\Lambda'}$ has trivial central character and so 
gives a representation of $\PU(\Phi)_\infty;$ and the above theorem is applicable to $\PU(\Phi)$. (One expects, although we have not carried this out, that Clozel's theorem should apply to reductive groups as well and so in particular to $\Unit(\Phi)$.) Having gotten a representation of $\PU(\Phi)$ we will think of it as a representation of 
$\Unit(\Phi)$ with trivial central character. We get the following consequence of Clozel's theorem: 

\smallskip

\begin{prop}
\label{prop:clozel}
Let $F, G$ and $\mu$ be as in Theorem~\ref{thm:main}.  
Fix two distinct finite places $v$ and $w$ of $F_0.$ Then there exists a cuspidal automorphic representation 
$\pi'$ of $G'/F_0$ such that
\begin{itemize}
\item $\pi'_\infty = \pi_{\Lambda'}$, the discrete series representation of $G'_\infty = G'(F_0 \otimes \R)$ as in 
Proposition~\ref{prop:cohomology},  
\item $\pi'_v$ is a supercuspidal representation of $G'(F_{0,v})$, and furthermore 
\item 
  \begin{itemize}
  \item in Cases (1a) and (1b): $\pi'_w$ is the Steinberg representation of $G'(F_{0,w}),$ whereas,  
  \item in Case (2): $v$ is chosen such that $\Unit(\Phi) \times F_{0,v} \simeq \GL(n)/F_{0,v},$ and we impose no condition on $w.$  
  \end{itemize}
\end{itemize}
In Case (2), the representation $\pi'$ has trivial central character. 
\end{prop}

\medskip
\subsubsection{\bf Consequences of the classification of the discrete spectrum of classical groups}

Let's begin by recalling the following result due to Magaard and Savin \cite[Proposition 8.2]{Magaard-Savin}: 

\begin{prop}
\label{prop:magaard-savin-sp2n}
Let $\sigma$ be a cuspidal automorphic representation on $\Sp(2n)$ over a number field $F$, 
such that for some finite place $v$ of $F$ the local component $\sigma_v$ is 
the Steinberg representation. Let $\pi$ be the automorphic representation of $\GL(2n+1)/F$ which is 
the lift of $\sigma$ as in \cite[Theorem 1.5.2]{Arthur}. Then $\pi_v$ is the Steinberg representation and 
$\pi$ is cuspidal.
\end{prop}

The key point in the proof of the above proposition in {\it loc.\,cit.,}
is to show that $\pi_v$ is the Steinberg representation (cuspidality of $\pi$ then 
easily follows). But, it is not {\it a priori} guaranteed from Arthur's work whether the Steinberg representation $\sigma_v$ transfers to the Steinberg representation. However, Arthur transfer is well-behaved with respect to the Aubert involution, and this involution switches the Steinberg representation with the trivial representation, and one observes using strong approximation for $\Sp(2n)$ that an automorphic representation containing the trivial representation as some local component is itself the trivial representation. Let's remark that this proof goes through {\it mutatis mutandis} for the split odd orthogonal group after one makes the same observation concerning the trivial representation for the split odd orthogonal group: let $\tau$ be an automorphic representation of $\SO(2n+1)$ with a trivial local component, say at $v$; now inflate $\tau$ to an automorphic representation $\tilde\tau$ of ${\rm Spin}(2n+1)$ which also has a trivial component at $v$; applying strong approximation which is known for the almost simple simply-connected group ${\rm Spin}(2n+1)$
(see, for example, Platonov--Rapinchuk \cite{Platonov-Rapinchuk}) we conclude that $\tilde\tau$ is trivial, {\it a fortiori},  
$\tau$ is trivial. We record the analogue of the above proposition for odd orthogonal group as:

\begin{prop}
\label{prop:magaard-savin-so}
Let $\sigma$ be a cuspidal automorphic representation of the split $\SO(2n+1)$ over a number field $F$, 
such that for some finite place $v$ of $F$ the local component $\sigma_v$ is 
the Steinberg representation. Let $\pi$ be the automorphic representation of $\GL(2n+1)/F$ which is 
the lift of $\sigma$ as in \cite[Theorem 1.5.2]{Arthur}. Then $\pi_v$ is the Steinberg representation and 
$\pi$ is cuspidal.
\end{prop}

\medskip
\subsubsection{\bf Conclusion of the proof}

Now consider the situation of Proposition~\ref{prop:clozel}. The cuspidal automorphic representation $\pi'$ of $G'$ transfers to an automorphic representation of $G = \GL_N/F$; this is via Arthur~\cite[Theorem 1.5.2]{Arthur} in cases (1a) and (1b) and by Mok\,\cite[Theorem 2.5.2]{Mok} in case (2). In cases (1a) and (1b), 
Propositions~\ref{prop:magaard-savin-sp2n} and \ref{prop:magaard-savin-so} ensure that $\pi$ is cuspidal. In case (2) it is clear since the place $v$ is taken so that the unitary group splits and has a supercuspidal local component guaranteeing cuspidality of $\pi.$ Now by Proposition~\ref{prop:cohomology} we know that $\pi$ is cohomological with respect to the given weight $\mu.$

\medskip

\section{Proof of Theorem~\ref{thm:auto-induction}}

We give a proof of Theorem ~\ref{thm:auto-induction} generalizing Clozel's construction in \cite{Clozel-Duke} using 
automorphic induction. 
Let $F$ be the given number field and $\mu$ be the given parallel weight assumed to have purity weight ${\sf w} = 0.$ 
So $\mu$ is of the form $\mu = (\mu_1 \geq \cdots \geq \mu_{2n})$ with $\mu_j \in \Z$ and  
such that 
$\mu_j  = -\mu_{2n+1-j}, \ \forall~ 1 \leq j \leq n$. Let $\rho = \left( n - \frac{1}{2}, n - \frac{3}{2}, \cdots, -n + \frac{1}{2}\right)$
 be the half sum of positive roots for $\mathfrak{gl}(2n,\C)$. Put $\ell = (\ell_1,\dots,\ell_{2n}) = 
 2\mu + 2\rho.$ Note that $\ell_{2n-i+1} = -\ell_i.$
  We proceed as in \cite{Clozel-Duke}. Choose a totally real number field $K_1$, which is 
 cyclic over $\Q$ and linearly disjoint with $F$ over $\Q$; let $\text{Gal}(K_{1}|\Q) = \left\langle \sigma \right\rangle \cong \Z / n \Z$. We also choose
a totally imaginary quadratic extension $L_1$ of $K_1$ which is also linearly disjoint with $F$ over $\Q.$ Let $L = FL_1$ and $K =FK_1$ be the corresponding compositum fields as shown in the diagram below. 

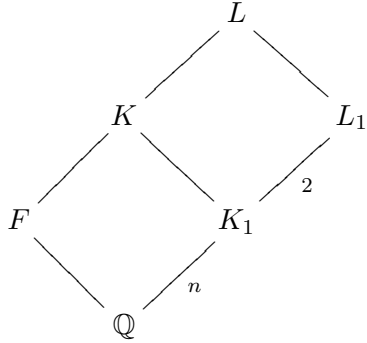
\hspace{2 cm}
\xymatrix{
& & & & & L \ar@{-}[dl] \ar@{-}[dr] \\
& & & & K \ar@{-}[dl] \ar@{-}[dr] & & L_{1} \ar@{-}[dl]^{2} \\
& & &F \ar@{-}[dr] &  & K_1 \ar@{-}[dl]^{n}\\
& & & &\Q }
\bigskip
 
Let $\I_{L_1} = L_1^\times\backslash \A_{L_1}^{\times}$ be the id\`ele class group of the number field $L_1$. 
Let $\chi_1 \in \text{Hom}(\I_{L_{1}},S^{1})$ be a unitary algebraic Hecke character of $\I_{L_{1}}$. From
 \cite{Weil}, one knows that the inifinity type of $\chi_1,$ i.e., its restriction
  to the group $\prod \limits_{v \in S_{\infty}(L_{1})} L_{1,v}^{\times}
  \cong (\C^{\times})^{n},$ is of the form: 
 \[ 
 \chi_{1}((a_v)) = \prod \limits_{v \in S_{\infty}(L_{1})} \left( a_{v}/\left|a_v \right| \right)^{-f_{v}},   
 \]  
 for integers $f_v.$ Conversely, from \cite{Weil}, given any set of integers $\{f_v\}$, there is a Hecke character with infinity type as above; this is because $L_1$ is a totally imaginary quadratic extension of a totally real
number field $K_1$. Fix an ordering on $S_{\infty}(L_{1})$ and take $\{f_v\} = (\ell_1,\dots,\ell_n),$ to get 
a character $\chi_1$ 
with infinity type: 
\[
(\chi_1) (z_1, z_2, \dots, z_n) =  \prod \limits_{j=1}^{n}( z_{j}/\overline{z_{j}})^{\ell_{j}/2} = 
\prod \limits_{j=1}^{n}( z_{j}/\left| z_{j} \right|)^{\ell_{j}}, 
\]
where, $(z_1,\dots,z_n) = (z_v)_{v \in S_{\infty}(L_{1})}.$ 
 
\medskip
Now we briefly explain the construction of a representation $\pi(\chi)$ of
$\text{GL}(2n,\A_{F})$ obtained by automorphic induction (in two steps) from $\chi =  \chi_{1}\circ N(L | L_{1})$. 
Following Jacquet--Langlands~\cite{Jacquet-Langlands}, one constructs a cuspidal
 representation $\pi_{K}(\chi)$ of $\text{GL}(2,\A_{K})$ via automorphic induction across $L/K.$ 
 Next, we let $\Pi_K$
 be the representation of $\text{GL}(2n,\A_{K})$ obtained by inducing
  $$
  \pi_{K}(\chi) \times \pi_{K}(\chi)^{\sigma} \times \cdots \times \pi_{K}(\chi)^{\sigma^{n-1}}
  $$
 from the Levi component $\text{GL}(2,\A_{K}) \times \text{GL}(2,\A_{K}) \times \cdots 
 \times \text{GL}(2,\A_{K})$ of the parabolic subgroup $\text{P}(2,2,\dots,2)$.
 Since $\Pi_K$ is stable under $\text{Gal}(K_{1}|\Q)$ action,
  it descends to a unique cuspidal representation 
 $\pi(\chi)$  of $\text{GL}(2n,\A_{F})$ such that for every archimedean place $v$ of $F$, 
 the local component
  $\pi_v$  is described as follows: 
  
 \begin{itemize}
 \item For a real place $v$ of $F$, 
       $$\pi_v = J_{v}(\mu): = \text{Ind}^{\rm{GL}(2n)}_{\text{P}(2,2,\cdots,2)} \left( D(\ell_1)        
         \otimes D(\ell_2) \otimes \cdots D(\ell_n) \right).$$
         
 \smallskip   
 \item For a complex place $v$ of $F$, 
       $$\pi_v = J_{v}(\mu): = \text{Ind}^{\rm{GL}_{2n}(\C)}_{\text{B}(\C)} \left( z^{a_{1}}\overline{z}^{b_{1}} \otimes
       \cdots \otimes  z^{a_{2n}}\overline{z}^{b_{2n}} \right). $$ 
       where $a: = \mu + \rho$, $b: =  - \mu - \rho$. 
    \end{itemize}
 
 As in \cite{Clozel-Duke}, from the results of Speh \cite{Speh} and Enright \cite{Enright}, 
 one knows that in both the above cases, $J_{v}(\mu)$ is the only
 \emph{generic} representation of $\GL(2n,F_v)$ which is cohomological with respect to the given parallel weight $\mu$. The proof of Theorem~\ref{thm:auto-induction} follows as in \cite{Clozel-Duke} 
 using the above observations and (\ref{eqn:cusp-coh}).

\medskip

\section{Remarks} 
\label{sec:remarks}

\subsection{Transferring from even orthogonal groups}
\label{sec:O(2n)toGL(2n)}
One may ask what sort of representations of $\GL_n$ are obtained by transferring from even orthogonal groups in as much as cohomological properties are concerned.  The basic problem is that this transfer does not preserve ``algebraicity." This is already seen in 
Gan--Raghuram \cite[Lemma 9.2\,(2)]{Gan-Raghuram} for local unramified representations. 
We begin by illustrating the issue if we use
Ramakrishnan's transfer \cite{Ramakrishnan-Annals} from $\GL(2) \times \GL(2)$ to $\GL(4)$,
which under a condition on central characters of the representations on $\GL(2) \times \GL(2),$ 
is the same as transferring from ${\rm SGO}(4)$ to $\GL_4$; see, for example, 
Ramakrishnan \cite[Section\,2]{Ramakrishnan-IMRN}.

\subsubsection{\bf On Ramakrishnan's transfer from $\GL(2) \times \GL(2)$ to $\GL(4)$}
For $j = 1,2$, let $\varphi_j$ be a holomorphic cuspidal normalized eigen-newform of weight $k_j.$ Let
$\pi(\varphi_j)$ be the associated cuspidal automorphic representation of $\GL_2(\A)$. Recall that
the representation at infinity is $\pi(\varphi_j)_\infty = D_{k_j-1}.$
Choose $\epsilon_j \in \{0,1\}$ such that $k_j \equiv \epsilon_j \pmod{2}$. Then
$\pi_j = \pi(\varphi_j) \otimes |\ |^{\epsilon_j/2}$ is a cohomological cuspidal representation.
(See \cite[Theorem\,5.5\,(2)]{Raghuram-Shahidi} and \cite[Theorem\,1.4]{Raghuram-Tanabe}.)
Consider the representation $\Pi = \pi_1 \boxtimes \pi_2$ as constructed in \cite{Ramakrishnan-Annals}. 
The representation at infinity is given by
$$
\Pi_\infty = D_{k_1-1}(\tfrac{\epsilon_1}{2}) \boxtimes D_{k_2-1}(\tfrac{\epsilon_2}{2}) \ = \
(D_{k_1-1} \boxtimes D_{k_2-1})(\tfrac{\epsilon_1 + \epsilon_2}{2}) \ = \
(D_{k_1+k_2-2} \times D_{k_1-k_2})(\tfrac{\epsilon_1 + \epsilon_2}{2}), 
$$
where the right-most representation is obtained by parabolic induction from the $(2,2)$-parabolic subgroup 
in $\GL_4(\R).$ Now suppose there is a dominant integral pure weight $\mu = (\mu_1, \mu_2, \mu_3, \mu_4)$, with
$\mu_1 + \mu_4 = \mu_2 + \mu_3 = {\sf w}$ such that $\Pi$ is cohomological with respect to $\M_\mu^{\sf v}$ then: 
$$
\Pi_\infty \ = \ D_{2\mu_1+3-{\sf w}}(\tfrac{{\sf w}}{2}) \times D_{2\mu_2+1-{\sf w}}(\tfrac{{\sf w}}{2}).
$$
This would entail the following equalities: 
${\sf w} \ = \ \epsilon_1 + \epsilon_2,$ $2\mu_1+3-{\sf w} = k_1+k_2-2,$ and $2\mu_2+1-{\sf w} = k_1-k_2,$
which is impossible for parity reasons. A possible remedy is suggested by Clozel \cite{Clozel-Motives} by
considering appropriate half-integral Tate twists. Define
$$
\Pi^T \ := \ \pi_1 \stackrel{T}{\boxtimes} \pi_2 \ := \ (\pi_1(-\tfrac{1}{2}) \boxtimes \pi_1(-\tfrac{1}{2}))(\tfrac{3}{2}) \ = \
(\pi_1 \boxtimes \pi_2)(\tfrac{1}{2}).
$$
Now, we would have
$$
\Pi^T_\infty \ = \ (D_{k_1+k_2-2} \times D_{k_1-k_2})(\tfrac{\epsilon_1 + \epsilon_2+1}{2})
\ = \ D_{2\mu_1+3-{\sf w}}(\tfrac{{\sf w}}{2}) \times D_{2\mu_2+1-{\sf w}}(\tfrac{{\sf w}}{2}).
$$
This is possible by taking
$
{\sf w}  = \epsilon_1 + \epsilon_2+1, 
\mu_1= (k_1+k_2+ \epsilon_1 + \epsilon_2)/2 - 2$ and 
$\mu_2 =  (k_1-k_2 + \epsilon_1 + \epsilon_2)/2.$
There are similar difficulties in the general case of transferring from even orthogonal groups for which there seems to be no obvious remedy.

\smallskip
\subsubsection{\bf Difficulties in the general case}
Let $N = 2n,$ and for simplicity take $F = \Q.$ Suppose that we want to construct a representation as in Case (1b) of the proof of Theorem \ref{thm:main} using an even orthogonal group as our endoscopy group. 
We follow \cite{Gross-Reeder}: take an even orthogonal group $G'/\Q$ such that 
if $n$ is even then $G'(\R) = \SO(n,n),$ and if $n$ is odd $G'(\R) = \SO(n-1,n+1).$ In both cases $G'(\R)$ has 
discrete series representations. Let's denote $G' =  \SO(2n)/\Q$ (resp., $G' =  \SO'(2n)/\Q$) if $n$ is even (resp., odd). 
Let $ \mu = (\mu_1 \geq \mu_2 \geq \cdots \geq \mu_n \geq -\mu_n \geq \cdots
\geq -\mu_2 \geq -\mu_1)$ be the given dominant integral weight for $\GL_N$ with purity weight $0$. Let $\mu'$, $\rho'$, $\Lambda'$ and $\ell'$ be defined as follows:
\begin{equation*}
\begin{array}{ll}
\mu'  & =  \ (\mu_1, \mu_2, \dots \mu_n);  \\ \smallskip
\rho' & = \ (n-1, n-2, \dots, 1, 0); \\ \smallskip
\Lambda'  & = \   \mu' + \rho' = (\mu_1+n - 1, \ \mu_2+n - 2, \dots, \ \mu_{n-1}+1, \ \mu_n); \\ \smallskip
\ell'  & = \ (\ell'_i)_{i=1}^{n}: \ = \ 2 \Lambda' \ = \  
(2\mu_1+2n - 2, \  2\mu_2 + 2n - 4,~ \dots, \  2\mu_{n-1}+2, \ 2\mu_n).
\end{array}
\end{equation*}
If we use a similar argument, we get a representation $\pi_\infty$ of $\GL_{2n}(\R)$ given by: 
$$
\pi_{\infty} = \Ind^{\GL_N(\R)}_{P(2,2,\dots, 2)} \left( D_{ \ell'_1} \otimes
 D_{ \ell'_2} \otimes \cdots \otimes  D_{ \ell'_n} \right), 
 $$
which is in fact not cohomological with respect to the weight $\mu.$ The {\it generic} representation of 
$\GL_{2n}(\R)$ that has nonzero cohomology with respect to $\mu^{\sf v}$ is given by inducing from 
$D_{ \ell_1}\otimes  D_{ \ell_2} \otimes \cdots \otimes  D_{ \ell_n},$ where 
$(\ell_1,\dots, \ell_{2n}) = 2\mu + 2\rho $ with $\rho$ being the half sum of positive roots for $\GL_N$. This gives  
 $(\ell_1,\dots, \ell_n) =  (2\mu_1 + 2n - 1, \ 2\mu_2 + 2n - 3,  \dots, 2\mu_n+1).$
Observe that $(\ell_1,\dots, \ell_n)$ and $(\ell_1',\dots, \ell_n')$ differ by $1$ at every coordinate. The remedy of taking a half-integral Tate twist in the previous subsection will not work here. It seems (to the authors) that there is no easy way to resolve this difficulty; the question being how one might modify the Langlands transfer from $\SO(2n)$ or $\SO'(2n)$ so as to preserve the property of being cohomological.

\medskip
\subsection{An endoscopic stratification of inner cohomology} 
Let's recall the definition of inner or interior cohomology. Take a field $E$ that is Galois over $\Q$ and containing a copy of $F.$ Refining the notations as in Section~\ref{sec:cuspidal-cohomology}, consider a dominant integral weight 
$\mu \in X^+(T \times E)$, giving us a rational finite-dimensional representation $\M_{\mu,E}$ of $G \times E$, which in turn gives a sheaf $\tM_{\mu,E}$ of $E$-vector spaces on $S^G_{K_f}.$ Inner cohomology is the image of cohomology with compact supports in global cohomology: 
$$
H^{\bullet}_!(S^{G_n}_{K_f}, \tM_{\mu,E}^{\sf v}) \ := \ 
{\rm Image}\left(H^{\bullet}_c(S^{G_n}_{K_f}, \tM_{\mu,E}^{\sf v})  \ \to \ 
H^{\bullet}(S^{G_n}_{K_f}, \tM_{\mu,E}^{\sf v}) \right). 
$$
If we pass to a transcendental situation via any embedding $\iota : E \to \C$, then it is well-known 
$$
H^{\bullet}_{\rm cusp}(S^{G_n}, \tM_{{}^\iota\!\mu,\C}^{\sf v}) \ \subset \ 
H^{\bullet}_!(S^{G_n}_{K_f}, \tM_{\mu,E}^{\sf v}) \otimes_{E,\iota} \C.  
$$
(See, for example, \cite{Clozel-Motives} or \cite{Harder-Raghuram}.) Next, let's recall arithmeticity for Shalika models (\cite[Appendix]{Grobner-Raghuram}): if  $\Pi$ is a cuspidal representation of $\GL_{2n}$ over a totally real $F$, and suppose that $\Pi$ is cohomological and has a Shalika model, then for $\sigma \in \autc$, the representation ${}^\sigma\Pi$ also has a Shalika model.  In other words, if $\Pi$ is cohomological and is a transfer from $\SO(2n+1)$ then so is any conjugate of $\Pi.$ The above considerations gives the following  

\begin{cor}
Let $F$ be a totally real field, and take $G = \GL(2n)/F.$ Let $\mu$ be a parallel weight with purity $0.$ Let $E$ be a Galois extension of $\Q$ that contains a copy of $F$. There exists a nontrivial $E$-subspace 
$$
H^\bullet_{\rm symp}(S^G, \M_{\mu, E}^{\sf v}) \ \subset \ H^\bullet_{!}(S^G, \M_{\mu, E}^{\sf v})
$$
stable under all Hecke operators, such that $H^\bullet_{\rm symp}(S^G, \M_{\mu, E}^{\sf v}) \otimes_{E, \iota} \C$ 
is spanned by cuspidal representations of $G$ that are all transfers from ${\rm SO}(2n+1).$
\end{cor}

Similarly, another result on arithmeticity (\cite[Remark 5.5]{Gan-Raghuram}), says that 
if a cuspidal representation $\Pi$ of $\GL(2n+1)$ over a totally real field 
is such that a partial symmetric square $L$-function has a pole at $s=1$ (a property that characterizes $\Pi$ being a transfer from $\Sp(2n)$) then the same is true of the representation ${}^\sigma\Pi$ for $\sigma \in \autc.$ This gives us

\begin{cor}
Let $F$ be a totally real field, and take $G = \GL(2n+1)/F.$ Let $\mu$ be a parallel weight with purity $0.$ Let $E$ be a Galois extension of $\Q$ that contains a copy of $F$. There exists a nontrivial $E$-subspace
$$
H^\bullet_{\rm orth}(S^G, \M_{\mu, E}^{\sf v}) \ \subset \ H^\bullet_{!}(S^G, \M_{\mu, E}^{\sf v})
$$
stable under all Hecke operators, 
such that $H^\bullet_{\rm orth}(S^G, \M_{\mu, E}^{\sf v}) \otimes_{E, \iota} \C$ is spanned by cuspidal representations of $G$ that are all transfers from $\Sp(2n).$
\end{cor}

The notation $H^\bullet_{\rm symp}$ (resp., $H^\bullet_{\rm orth}$) is to suggest that we are looking at the contribution of cuspidal representations with `symplectic' (resp., `orthogonal') parameters. 

\bigskip

{\it Acknowledgements:} The authors thank Laurent Clozel, Wee Teck Gan, Dipendra Prasad and Gordan Savin for helpful conversations which have found their way into this article.

\bigskip

\bigskip

\end{document}